\documentclass[a4paper,12pt,leqno]{amsart}
\usepackage{amsmath,amsthm,amssymb,latexsym,epsfig,graphicx,subfigure}
\numberwithin{equation}{section}

\newtheorem{thm}{Theorem}[section]

\theoremstyle{definition}

\theoremstyle{definition}
\newtheorem{rem}[thm]{Remark}
\newtheorem{notation}[thm]{Notation}

\newcommand{\calH}{\mathcal{H}}

\newcommand{\R}{\mathbb{R}}

\newcommand{\N}{\mathbb{N}}

\newcommand {\grtrsim} {\ {\raise-.5ex\hbox{$\buildrel>\over\sim$}}\ }

\newcommand{\khii}{\text{\lower -.4ex\hbox{$\chi$}}}
\DeclareMathOperator{\spt}{spt}

\begin{document}
\title [Exceptional sets for visibility]{A sharp exceptional set estimate for visibility}

\thanks{TO is supported by the Academy of Finland through the grant \emph{Restricted families of projections and connections to Kakeya type problems.}} 
\author{Tuomas Orponen}

\begin{abstract} A Borel set $B \subset \R^{n}$ is \emph{visible} from $x \in \R^{n}$, if the radial projection of $B$ with base point $x$ has positive $\calH^{n - 1}$ measure. I prove that if $\dim B > n - 1$, then $B$ is visible from every point $x \in \R^{n} \setminus E$, where $E$ is an exceptional set with dimension $\dim E \leq 2(n - 1) - \dim B$. This is the sharp bound for all $n \geq 2$.

Many parts of the proof were already contained in a recent previous paper by P. Mattila and the author, where a weaker bound for $\dim E$ was derived as a corollary from a certain slicing theorem. Here, no improvement to the slicing result is obtained; in brief, the main observation of the present paper is that the proof method gives the optimal result, when applied directly to the visibility problem. 
\end{abstract}

\maketitle

\section{Introduction}

For $x \in \R^{n}$, let $\pi_{x} \colon \R^{n} \setminus \{x\} \to S^{n - 1}$ be the radial projection
\begin{displaymath} \pi_{x}(y) := \frac{y - x}{|y - x|}. \end{displaymath}
The following theorem is the main result of the paper:
\begin{thm}\label{main} Assume that $B \subset \R^{n}$ is a Borel set with $\dim B > n - 1$. Then, there exists a set $E \subset \R^{n}$ with $\dim E \leq 2(n - 1) - \dim B$ such that 
\begin{displaymath} \calH^{n - 1}(\pi_{x}(B)) > 0, \qquad x \in \R^{n} \setminus E. \end{displaymath}
This is the sharp bound for every $n \geq 2$. 
\end{thm}
This settles a conjecture by P. Mattila and the author in \cite{MO}, where it was proven that $\dim E \leq n - 1$ as soon as $\dim B > n - 1$. The same conjecture had earlier appeared in Mattila's survey paper \cite{Mat1}, see (6.1) on p. 36.\footnote{In fact, (6.1) in \cite{Mat1} states the planar version of Theorem \ref{main} not as a conjecture, but a \emph{fact}, which should follow from Peres and Schlag's work \cite{PS}. However, Theorem 7.3 in \cite{PS} gives the bound $3 - \dim B$ instead of $2 - \dim B$ in the plane.} Finally, the proof in the present paper also fills a small gap in the argument in \cite{MO}, see the footnote on page 3.
\begin{rem} The strict inequality $\dim B > n - 1$ is necessary. In fact, if $B \subset \R^{n}$ is purely $(n - 1)$-unrectifiable with $0 < \calH^{n - 1}(B) < \infty$, then it follows easily from the Besicovitch-Federer projection theorem that $\calH^{n - 1}(\pi_{x}(B)) = 0$ for almost every $x \in \R^{n}$. A more precise result is due to Marstrand \cite{Mar}, Theorem VI: if $B \subset \R^{2}$ is purely $1$-unrectifiable with $0 < \calH^{1}(B) < \infty$, then $\calH^{1}(\pi_{x}(B)) = 0$ for all $x \in \R^{2} \setminus E$, where $\dim E \leq 1$. It is entertaining to note that when $\dim B > 1$, the same is true with "$\calH^{1}(\pi_{x}(B)) = 0$" replaced by "$\calH^{1}(\pi_{x}(B)) > 0$". 

Marstrand's result can be further improved for self-similar sets: Simon and Solomyak \cite{SS} have shown that if $B \subset \R^{2}$ is a purely $1$-unrectifiable self-similar set in the plane with $0 < \calH^{1}(B) < \infty$, and satisfying the open set condition, then $\calH^{1}(\pi_{x}(B)) = 0$ for every base point $x \in \R^{2}$. There is also a recent, more quantitative, version of this result by Bond, \L{}aba and Zahl \cite{BLZ}. \end{rem}

\begin{notation} The Grassmannian manifold of all $(n - 1)$-dimensional subspaces of $\R^{n}$ is denoted by $G(n,n - 1)$, and the Haar probability measure on $G(n,n - 1)$ is denoted by $\gamma_{n,n - 1}$. Given a plane $V \in G(n,n - 1)$, the mapping $\pi_{V} \colon \R^{n} \to V$ is the orthogonal projection onto $V$. If $\mu$ is a Radon measure on $\R^{n}$, its push-forward under $\pi_{V}$ is denoted by $\pi_{V\sharp}\mu$, so that
\begin{displaymath} \pi_{V\sharp}\mu(B) = \mu(\pi_{V}^{-1}(B)), \qquad B \subset V. \end{displaymath}
For $a,b > 0$, we write $a \lesssim b$, if there exists a constant $C \geq 1$ such that $a \leq Cb$; the constant $C$ may, without special mention, depend on various "fixed" parameters in the proof, such as the dimension of the ambient space, or that of $B$ in Theorem \ref{main}. 

For $0 \leq s \leq n$, the $s$-dimensional Hausdorff measure is denoted by $\calH^{s}$. The notation $\dim$ stands for Hausdorff dimension. Finally, if $\mu$ is a Radon measure on $\R^{n}$ and $0 \leq s \leq n$, the $s$-energy of $\mu$ is denoted by $I_{s}(\mu)$, so that by definition
\begin{displaymath} I_{s}(\mu) = \iint \frac{d\mu x \, d\mu y}{|x - y|^{s}}. \end{displaymath}
It is well-known that, see Theorem 3.10 in \cite{Mat2}, that
\begin{equation}\label{form8} I_{s}(\mu) = c(n,s) \int |\hat{\mu}(\xi)|^{2}|\xi|^{s - n} \, d\xi, \qquad 0 < s < n. \end{equation}
\end{notation}

\section{Acknowledgements} This paper would not exist without the previous joint work of Pertti Mattila and myself in \cite{MO}, and I am grateful for his collaboration and insights. In addition, I thank Pertti for reading the current manuscript, and giving plenty of good comments.

\section{Proof of the main theorem} 

The first part of this section contains the proof of the bound $\dim E \leq 2(n - 1) - \dim B$. The second, far shorter, part discusses the question of sharpness.

\subsection{Proof of the upper bound for $\dim E$} It suffices to prove the theorem for compact sets $B$, because if $\dim \{x : \calH^{n - 1}(\pi_{x}(B)) = 0\} > 2(n - 1) - \dim B$ for some Borel set $B$, then also $\dim \{x : \calH^{n - 1}(\pi_{x}(K)) = 0\} > 2(n - 1) - \dim K$ for some compact set $K \subset B$ with $n - 1 < \dim K \leq \dim B$. So assume that $B$ is compact. Then, the set $E := \{x \in \R^{n} : \calH^{n - 1}(\pi_{x}(B)) = 0\}$ is Borel, and we make the counter assumption that
\begin{displaymath} 2(n - 1) - s < \dim E < n - 1 \end{displaymath}
for some $n - 1 < s < \dim B$ (such an $s$ can be found if $\dim E > 2(n - 1) - \dim B$). We may further assume that $E$ and $B$ are disjoint; if this is not true to begin with, choose two disjoint closed balls $B_{1}$ and $B_{2}$ such that $\dim [B \cap B_{1}] > s$, and $2(n - 1) - \dim [B \cap B_{1}] < \dim [E \cap B_{2}]$. Finally, fix $t$ strictly between $2(n - 1) - s$ and $\dim E$, and find compactly supported Borel probability measures $\mu$ and $\nu$ inside $B$ and $E$, respectively, such that $I_{s}(\mu) < \infty$, and $I_{t}(\nu) < \infty$. Then $\calH^{n - 1}(\pi_{x}(\spt \mu)) = 0$ for every $x \in \spt \nu$; to simplify notation, we assume that $B = \spt \mu$ and $E = \spt \nu$.

We briefly discuss the meaning of the assumption $\calH^{n - 1}(\pi_{x}(B)) = 0$ for $x \in E$. If $L_{V,x}$ is the line perpendicular to $V \in G(n,n - 1)$ and passing through $x \in \R^{n}$, another way to write $\calH^{n - 1}(\pi_{x}(B)) = 0$, $x \in E$, is the following:
\begin{equation}\label{form7} \gamma_{n,n - 1}(\{V : L_{V,x} \cap B \neq \emptyset\}) = 0, \qquad x \in E. \end{equation}
This is where we needed to know that $B$ and $E$ are disjoint. Using Fubini's theorem, \eqref{form7} implies that
\begin{equation}\label{form1} \nu(\{x : L_{V,x} \cap B \neq \emptyset\}) = 0 \end{equation}
for $\gamma_{n,n - 1}$ almost every $V \in G(n,n - 1)$. 

For $\delta \in (0,1)$, let $\psi_{\delta} \colon \R^{n} \to [0,\infty)$ be a radial compactly supported approximate identity (thus $\psi_{\delta} = \delta^{-n}\psi(x/\delta)$, where $\psi$ is non-negative, radial, supported on $B(0,1)$ and has integral one). Let $\mu_{\delta} := \mu \ast \psi_{\delta}$, and consider the function 
\begin{displaymath} V \mapsto f_{\delta}(V) := \int_{V} \pi_{V\sharp} \mu_{\delta} \, d\pi_{V\sharp}\nu, \qquad V \in G(n,n - 1). \end{displaymath}
We will need to know that
\begin{itemize}
\item[(i)] $\|f_{\delta}\|_{L^{1}(G(n,n - 1))} \geq c$ for some constant $c > 0$ independent of $\delta \in (0,1)$, and
\item[(ii)] there exists $p > 1$ (depending on $n,s$ and $t$ only) such that $\|f_{\delta}\|_{L^{p}(G(n,n - 1))} \leq C$, where $C < \infty$ is independent of $\delta \in (0,1)$.\footnote{This $L^{p}$-estimate was missing from the paper \cite{MO}.}
\end{itemize}
In fact, (i) is precisely (3.4) in \cite{MO}, so we skip the details: in brief, applying the Parseval formula and integrating in polar coordinates, one can show that $\|f_{\delta}\|_{L^{1}}$ equals a constant times $\iint |x - y|^{-(n - 1)} \, d\mu_{\delta}x \, d\nu x$, which is uniformly bounded from below for $\delta \in (0,1)$.

We then prove (ii). Write $s' := 2(n - 1) - t < s$, and let $\sigma$ be a measure on $G(n,n - 1)$ satisfying the growth condition $\sigma(B(x,r)) \lesssim r^{h}$ for some 
\begin{displaymath} \max\{t,n - 1 + s' - s\} < h < n - 1. \end{displaymath}
Write $\mu_{V}^{\delta} := \pi_{V \sharp}\mu_{\delta}$, and $\nu_{V} := \pi_{V\sharp}\nu$. Under the previous restrictions, it is known (see discussion below) that
\begin{equation}\label{form3} \int \int_{V} |x|^{t - (n - 1)}|\widehat{\nu_{V}}(x)|^{2} \, d \calH^{n - 1}(x)\, d\sigma V = c(n,t)\int I_{t}(\nu_{V}) \, d\sigma V \lesssim I_{t}(\nu) < \infty
\end{equation}
and
\begin{equation}\label{form5}\int \int_{V} |x|^{s' - (n - 1)}|\widehat{\mu^{\delta}_{V}}(x)|^{2} \, d\calH^{n - 1}(x) \, d\sigma V \lesssim 1 + I_{s}(\mu_{\delta}) \lesssim 1 + I_{s}(\mu) < \infty. \end{equation}
The bound \eqref{form3} in the plane is due to Kaufman \cite{Kau}, and the higher dimensional analogue we need can be found in a paper by Mattila, see Lemma 5.1 in \cite{Mat}. As such, the bound \eqref{form5} is most likely due to Peres and Schlag \cite{PS}, but it is certainly inspired by earlier work of Falconer \cite{Fal1}; a proof can also be found on p. 81 in the book \cite{Mat2}.

Armed with Parseval, \eqref{form3}, \eqref{form5} and Cauchy-Schwarz, we make the following estimate:
\begin{align} \int & \int_{V} \mu_{V}^{\delta} \, d\nu_{V} \, d\sigma V \leq \int \int_{V} |\widehat{\mu^{\delta}_{V}}(x)||\widehat{\nu_{V}}(x)| \, d\calH^{n - 1}(x) \, d\sigma V\notag \\
&  = \int \int_{V} |x|^{(s' - (n - 1))/2}|x|^{(t - (n - 1))/2}|\widehat{\mu^{\delta}_{V}}(x)||\widehat{\nu_{V}}(x)| \, d\calH^{n - 1}(x) \, d\sigma V \notag\\
& \label{form4}\leq \int \left(\int_{V} |x|^{t - (n - 1)}|\widehat{\nu_{V}}(x)|^{2}\,d\calH^{n - 1}(x) \right)^{1/2} \left(\int_{V} |x|^{s' - (n - 1)}|\widehat{\mu^{\delta}_{V}}(x)|^{2}\,d\calH^{n - 1}(x) \right)^{1/2} \, d\sigma V\\
& \leq \left(\int \int_{V} |x|^{t - (n - 1)}|\widehat{\nu_{V}}(x)|^{2}\,d\calH^{n - 1}(x) \, d\sigma V \right)^{1/2}\notag \\
&\qquad \times \left( \int \int_{V} |x|^{s' - (n - 1)} |\widehat{\mu^{\delta}_{V}}(x)|^{2}\,d\calH^{n - 1}(x) \, d\sigma V \right)^{1/2} \lesssim I_{t}(\nu)^{1/2}(1 + I_{s}(\mu))^{1/2}.\notag \end{align} 
Next, to get the $L^{p}$-result we desired, we observe that functions $g \in L^{q}(G(n,n - 1))$ with $\|g\|_{L^{q}} = 1$ (where $q$ is dual to $p$) satisfy the kind of "power bound" as was required of $\sigma$. Namely,
\begin{displaymath} \int_{B(V,r)} g \, d\gamma_{n,n - 1} \leq \gamma_{n,n - 1}(B(V,r))^{1/p} \left( \int |g|^{q} \, d\gamma_{n,n - 1} \right)^{1/q} \lesssim r^{(n - 1)/p}. \end{displaymath}
So, if $p > 1$ is so close to one that that $(n - 1)/p \geq h$, the estimate \eqref{form4} yields
\begin{align*} \int f_{\delta} \cdot g \, d\gamma_{n,n - 1} & \leq \int \left( \int_{V} |\widehat{\mu^{\delta}_{V}}(x)||\widehat{\nu_{V}}(x)| \, d\calH^{n - 1}(x) \right) g(V) \, d\gamma_{n,n - 1}(V)\\
& \lesssim I_{t}(\nu)^{1/2}(1 + I_{s}(\mu))^{1/2}. \end{align*}
By the usual $L^{p}-L^{q}$ duality, this proves (ii). 

We record two further standard facts: for $\gamma_{n,n - 1}$ almost every $V \in G(n,n - 1)$,
\begin{itemize}
\item[(iii)] the measure $\pi_{V\sharp}\mu$ lies in the fractional Sobolev space $H^{(s - (n - 1))/2}(V)$, and
\item[(iv)] the measure $\pi_{V\sharp}\nu$ has finite $t$-energy.
\end{itemize}
Fact (iv) follows immediately from \eqref{form3} with $\sigma = \gamma_{n,n - 1}$. Fact (iii) does not quite follow from \eqref{form5} as stated above (because we assumed $s' < s$), but it does follow from the variant of \eqref{form5}, where $s' = s$ and $\sigma = \gamma_{n,n - 1}$; this remains true, as can be proven easily via "integration in polar coordinates", see for instance (24.2) in \cite{Mat2}. This gives fact (iii). 

Assume that $V \in G(n,n - 1)$ is a plane such that (iii) and (iv) hold. Then, as observed already in \cite{MO} (or see Theorem 17.3 in the book \cite{Mat2}), the Hardy-Littlewood maximal function of $\pi_{V\sharp}\mu$ belongs to $L^{1}(\pi_{V\sharp}\nu)$, which implies that the functions $\pi_{V\sharp}\mu_{\delta}$ converge to a limit in $L^{1}(\pi_{V\sharp}\nu)$ both in $L^{1}(\pi_{V\sharp}\nu)$, and $\pi_{V\sharp}\nu$ almost everywhere (to see this, observe that $\pi_{V\sharp}\mu_{\delta} = \psi^{V}_{\delta} \ast \pi_{V\sharp}\mu$ for some approximate identity $\psi^{V}_{\delta}$ on $V$, because $\psi$ was chosen radial). We denote the limit by $g_{V} \in L^{1}(\pi_{V\sharp}\nu)$, so that
\begin{displaymath} [\pi_{V\sharp}\mu_{\delta}](v) \to g_{V}(v) \text{ for } \pi_{V\sharp}\nu \text{ almost every } v \in V. \end{displaymath}

Recalling from (ii) that the sequence $(f_{\delta})_{\delta > 0}$ is bounded in $L^{p}(G(n,n - 1))$ for some $p > 1$, we may pick a subsequence $(f_{\delta_{j}})_{j \in \N}$, which converges weakly in $L^{p}(G(n,n - 1))$ to a limit $f \in L^{p}(G(n,n - 1))$. The values of $f$ are known to us: since $\pi_{V\sharp}\mu_{\delta} \to g_{V}$ in $L^{1}(\pi_{V\sharp}\nu)$, whenever (iii) and (iv) hold (that is, for $\gamma_{n,n - 1}$ almost every $V$), the whole sequence $(f_{\delta})_{\delta > 0}$ converges pointwise $\gamma_{n,n - 1}$ almost everywhere, and we may infer that
\begin{displaymath} f(V) = \int_{V} g_{V} \, d\pi_{\sharp}\nu \text{ for } \gamma_{n,n - 1} \text{ almost every } V \in G(n,n - 1). \end{displaymath} 
On the other hand, since $f$ is the weak $L^{p}$-limit of the functions $f_{\delta_{j}}$, each of which satisfies the uniform $L^{1}$ lower bound from (i), we have 
\begin{displaymath} \|f\|_{L^{1}(G(n,n - 1))} \geq \limsup_{j \to \infty} \|f_{\delta_{j}}\|_{L^{1}(G(n,n - 1))} \geq c > 0. \end{displaymath}
This estimate is legitimate, because $G(n,n - 1)$ is compact. It follows that $f(V) > 0$ for $\gamma_{n,n - 1}$ positively many planes $V \in G(n,n - 1)$. 

Using this fact, we find a plane $V \in G(n,n - 1)$ with the following four properties: \eqref{form1}, (iii) and (iv) hold, and
\begin{equation}\label{form2} \int g_{V} \, d\pi_{V\sharp}\nu > 0. \end{equation}
The proof is finished by showing that the four conditions cannot, in fact, hold simultaneously. To this end, write
\begin{displaymath} G_{V} := \{v \in V : [\pi_{V\sharp}\mu_{\delta}](v) \to g_{V}(v)\}, \end{displaymath}
so that $\pi_{V\sharp}\nu(V \setminus G_{V}) = 0$, as we discussed after (iii) and (iv). Then, decompose the integral in \eqref{form2} as follows:
\begin{displaymath} \int g_{V} \, d\pi_{V\sharp}\nu = \int_{\{v \in G_{V} : \pi_{V}^{-1}\{v\} \cap B = \emptyset\}} \ldots + \int_{\{v \in G_{V} : \pi_{V}^{-1}\{v\} \cap B \neq \emptyset\}} \ldots + \int_{V \setminus G_{V}} \ldots. \end{displaymath}
The third integral is clearly zero, and the same is true for the second integral by \eqref{form1}:
\begin{displaymath} \pi_{V\sharp}\nu(\{v : \pi_{V}^{-1}\{v\} \cap B \neq \emptyset\}) = \nu(\{x : L_{V,x} \cap B \neq \emptyset\}) = 0 \end{displaymath}
But also the first integral is zero: indeed, if $v \in V$ and $\pi_{V}^{-1}\{v\} \cap B = \emptyset$, then the compactness of $B$ implies that $\pi_{V}^{-1}(B(v,\delta)) \cap B = \emptyset$ for $\delta > 0$ small enough. If moreover $v \in G_{V}$, this implies that
\begin{displaymath} g_{V}(v) = \lim_{\delta \to 0} [\pi_{V\sharp}\mu_{\delta}](v) = 0. \end{displaymath}
In other words, $g_{V}(v) = 0$ for every $v \in \{G_{V} : \pi_{V}^{-1}\{v\} \cap B = \emptyset\}$. We have now seen that \eqref{form2} cannot hold, and the ensuing contradiction completes the proof. 

\subsection{Sharpness of the bound} Given $n - 1 < s < n$, there exists a compact set $K \subset \R^{n}$ such that $\dim K = s$, and
\begin{equation}\label{form6} \dim \{V \in G(n,n - 1) : \calH^{n - 1}(\pi_{V}(K)) = 0\} = 2(n - 1) - s. \end{equation}
The example is due to Peltom\"aki \cite{Pe}, but the details can also be found in \cite{Mat2}, Example 5.13. 

Consider the projective transformation $F \colon \R^{n} \setminus \R^{n - 1} \to \R^{n}$, defined by
\begin{displaymath} F(\tilde{x},x_{n}) := \frac{(\tilde{x},1)}{x_{n}}, \qquad (\tilde{x},x_{n}) \in \R^{n - 1} \times \R \setminus \{0\}. \end{displaymath}
Then $F$ maps lines in $\R^{n}$ of the form $\{te + (a,0) : t \in \R\}$, where $a \in \R^{n - 1}$ and $e = (\tilde{e},e_{n}) \in S^{n - 1}$, $e_{n} \neq 0$, to lines of the form $\{u(a,1) + (\tilde{e}/e_{n},0) : u \in \R\}$. In particular, fixing the "base point" $a \in \R^{n - 1}$, the mapping $F$ takes the lines passing through $(a,0)$ to lines parallel to the vector $(a,1)$. Now, let $G := F^{-1}$, and consider the set $G(K) \subset \R^{n}$, which clearly still has $\dim G(K) = s$. The equation \eqref{form6} can be (essentially) reworded by saying that there exists a $2(n - 1) - s$ dimensional family $E$ of vectors of the form $(a,1)$ such that $K$ lies on an $\calH^{n - 1}$-null set of lines parallel to each $(a,1) \in E$. Hence, there exists a $2(n - 1) - s$ dimensional family $E'$ of points $a \in \R^{n - 1}$ such that $G(K)$ lies on an $\calH^{n - 1}$-null set of lines passing through $(a,0)$. In other words, $\pi_{(a,0)}(G(K)) = 0$ for every $a \in E'$, as desired.


\begin{thebibliography}{CMM}

\bibitem{Fal1}
K.J. Falconer. Hausdorff dimension and the exceptional set of projections, 
{\it Mathematika} {\bf 29} (1982), 109--115.

\bibitem{Kau}
R. Kaufman.
On Hausdorff dimension of projections,
 {\em Mathematika} {\bf 15} (1968), 153-155.
 
\bibitem{BLZ}
M. Bond, I. \L{}aba and J. Zahl.
Quantitative visibility estimates for unrectifiable sets in the plane,
appeared electronically in {\em Trans Amer. Math. Soc.} (2015)
 
 \bibitem{Mar}
 J.M. Marstrand. 
 Some fundamental geometrical properties of plane sets of fractional dimensions,
 {\em Proc. London Math. Soc.} (3) 4 (1954), 257--302

\bibitem{Mat}
P.~Mattila.
 Hausdorff dimension, orthogonal projections and intersections with
  planes,
 {\em Ann. Acad. Sci. Fenn. A Math.} {\bf 1} (1975),  227--244.
 
 \bibitem{Mat1}
 P. Mattila.
 Hausdorff dimension, projections, and the Fourier transform,
 {\em Publ. Mat.} \textbf{48} (2004), 3--48

\bibitem{Mat2}
P. Mattila.
Fourier Analysis and Hausdorff Dimension,
Cambridge University Press, Cambridge, 2015.

\bibitem{MO} P. Mattila and T. Orponen.
Hausdorff dimension, intersections of projections and exceptional plane sections, appeared electronically in \emph{Proc. Amer. Math. Soc.} (2015), available at arXiv:1509.05724

\bibitem{Pe} A. Peltom\"aki.
Projektiot ja Hausdorffin dimensio, Licenciate thesis, Helsingin yliopisto (1988)

\bibitem{PS} Y. Peres and W. Schlag. Smoothness of projections, Bernoulli convolutions, and the dimension 
of exceptions, {\em Duke Math. J.} {\bf 102} (2000), 193--251.   

\bibitem{SS} K. Simon and B. Solomyak. Visibility for self-similar sets of dimension one in the plane, {\em Real Anal. Exchange} {\bf 32}(1) (2006/2007), 67--78

\end{thebibliography}
\end{document}